\documentclass[leqno,12pt]{article}
\usepackage{amsmath,amssymb}
\usepackage{color}
\usepackage{epsfig}
\usepackage[utf8]{inputenc}
\usepackage[T1]{fontenc}

\newtheorem{Theorem}{\hspace{\parindent}\bf Theorem}[section]
\newtheorem{Lemma}{\hspace{\parindent}\bf Lemma}[section]
\newtheorem{Proposition}{\hspace{\parindent}\bf Proposition}[section]
\newtheorem{Corollary}{\hspace{\parindent}\bf Corollary}[section]

\newcommand{\qed}{\hfill$\square$\vspace{0.3cm}}

\newcommand{\R}{\mathbb{R}}
\newcommand{\ren}{\mathbb{R}^N}

\newcommand{\E}{\mathop{\rm E}}

\begin{document}

\title{\textbf{Classical solutions and higher regularity  for nonlinear fractional diffusion equations}}
\author{ by \\
Juan Luis V\'{a}zquez, Arturo de Pablo,  \\ Fernando Quir\'{o}s, and Ana Rodr\'{\i}guez}

\maketitle

\begin{abstract}
We study the regularity properties of the solutions to the nonlinear equation with fractional diffusion
$$
\partial_tu+(-\Delta)^{\sigma/2}\varphi(u)=0,
$$
posed for $x\in \mathbb{R}^N$, $t>0$, with $0<\sigma<2$, $N\ge1$. If the
nonlinearity  satisfies some not very restrictive
conditions: $\varphi\in C^{1,\gamma}(\mathbb{R})$, $1+\gamma>\sigma$, and $\varphi'(u)>0$ for every $u\in\mathbb{R}$, we prove that  bounded weak solutions are classical solutions for all positive times.  We also
explore sufficient conditions on the non-linearity to obtain higher
regularity  for the solutions, even $C^\infty$ regularity. Degenerate and singular cases, including the power nonlinearity $\varphi(u)=|u|^{m-1}u$, $m>0$, are also considered, and
the existence of classical solutions in the power case is proved.
\end{abstract}

\vskip 1cm

\noindent{\makebox[1in]\hrulefill}\newline

\newpage
\section{Introduction}\label{sect-introduction}
\setcounter{equation}{0}

This paper is devoted to establish the regularity of bounded weak
solutions for the nonlinear parabolic equation  involving fractional
diffusion
\begin{equation}  \label{eq:main}
\partial_tu+(-\Delta)^{\sigma/2}\varphi(u)=0  \quad\mbox{in } Q=\mathbb{R}^N\times(0,\infty).
\end{equation}
Here $(-\Delta)^{\sigma/2}=\mathcal{F}^{-1}(|\cdot|^\sigma\mathcal{F})$, where $\mathcal{F}$ denotes Fourier transform, is the usual fractional Laplacian with $0<\sigma<2$
and $N\ge1$. The constitutive function $\varphi$ is assumed to be at
least continuous and nondecreasing. Further
conditions will be introduced as needed.

The existence of a unique weak solution to the Cauchy problem for equation~\eqref{eq:main}
has been fully investigated in~\cite{pqrv, pqrv2} for the case where
$\varphi$ is a positive power. The solution in that case is in fact bounded
for positive times even if the initial data are not, provided they are in a suitable integrability space. Such
theory can be easily extended to the case of more general functions
$\varphi$; see Section~\ref{sec.exist} at the end of the paper for some details.

When $\varphi(u)=u$ the equation is the so-called fractional heat equation,
that has been studied in a number of papers, mainly coming from
probability. Explicit representation with a kernel allows to show in
this case that solutions are $C^\infty$ smooth and bounded for every
$t>0$, $x\in\mathbb{R}^N$, under the assumption that the initial data are integrable.
In the  nonlinear case such a representation is not available. Nevertheless, we will still be able to obtain that bounded weak solutions are smooth if the
equation is \lq\lq uniformly parabolic'',
$0<c\le \varphi'(u)\le C<\infty$.

\medskip

\noindent {\sc Classical solutions.} Our first result establishes that if the nonlinearity is smooth enough, compared to the order of the equation, $\max\{1,\sigma\}$, then bounded weak solutions are classical solutions.

\begin{Theorem}\label{th:main} Let $u$ be a bounded weak solution to \eqref{eq:main}, and assume $\varphi\in C^{1,\gamma}(\mathbb{R})$,
$0<\gamma<1$, and $\varphi'(s)>0$ for every
$s\in\mathbb{R}$. If $1+\gamma>\sigma$, then
$\partial_tu$ and $(-\Delta)^{\sigma/2}\varphi(u)$ are H\"older continuous functions and~\eqref{eq:main} is satisfied everywhere.
\end{Theorem}
The precise regularity of the solution is determined by the regularity of the nonlinearity $\varphi$; see Section~\ref{sec.beyond.Holder}
for the details. Notice that the condition $\varphi'>0$ together with the boundedness of $u$ implies that the equation is uniformly parabolic.

The idea of the proof is as follows: thanks to the results of Athanasopoulos and Caffarelli~\cite{Athanasopoulos-Caffarelli}, we already know that bounded weak solutions
are $C^\alpha$ regular for some $\alpha\in(0,1)$. In order to improve this regularity we write the  equation~\eqref{eq:main} as a fractional linear heat equation with a source term. This term is in principle not very smooth, but it has some good properties. To be precise, given  $(x_0,t_0)\in Q$, we have
\begin{equation}\label{eq:nonlin.lin}
\partial_tu+(-\Delta)^{\sigma/2} u=(-\Delta)^{\sigma/2} f,
\end{equation}
where
\begin{equation*}
\label{eq:linear-f}
f(x,t):=u(x,t)-\frac{\varphi(u(x,t))}{\varphi'(u(x_0,t_0))},
\end{equation*}
after the time rescaling $t\to t/\varphi'(u(x_0,t_0))$.
It turns out,  as we will prove in Sections~\ref{sec.regularity.nonlinear} and \ref{sec.beyond.Holder}, that  solutions to the linear equation~\eqref{eq:nonlin.lin} have the same regularity as $f$.

Next, using the nonlinearity we observe that  $f$ in the actual right-hand side is more regular than $u$ near $(x_0,t_0)$; see formula~\eqref{eq:local.holder.condition0}. We are thus in a situation that is somewhat similar to the one considered by Caffarelli and Vasseur in \cite{Caffarelli-Vasseur},  where they deal with an equation, motivated by the study of geostrophic equations, of the form
\begin{equation*}
\label{eq:caffa-vasseur}
\partial_tu+(-\Delta)^{1/2} u=\mathop{\rm div}( {\bf v}u),
\end{equation*}
where $\bf v$ is a divergence free vector. Comparing with \eqref{eq:nonlin.lin}, we see two differences: in their case $\sigma=1$, and the source term is local. These two differences will significantly complicate our analysis.

In order to obtain the above-mentioned regularity for the solutions $u$ to~\eqref{eq:nonlin.lin}, we will use the fact that they are given by the representation formula
\begin{equation*}
\label{eq:duhamel}
\begin{array}{l}
u(x,t)=\displaystyle\int_{\mathbb{R}^N}P_\sigma(x-\overline x, t) u(\overline x,0)\,d\overline x\\[4pt]
\displaystyle\qquad\qquad+\int_0^t\int_{\mathbb{R}^N}(-\Delta)^{\sigma/2}P_\sigma(x-\overline x,
t-\overline t) f(\overline x,\overline t)\,d\overline xd\overline t,\quad (x,t)\in Q,
\end{array}
\end{equation*}
where $P_\sigma$ is the kernel of the $\sigma$-fractional linear heat equation; see Section~\ref{sec.linear} for a proof of this fact, that falls into the linear theory. Therefore, we are led to study the singular kernel $A_\sigma(x,t):=(-\Delta)^{\sigma/2}P_\sigma(x,t)$. Unfortunately, $P_\sigma$, and hence $A_\sigma$, is only explicit when $\sigma=1$. However, using the self-similar structure of $P_\sigma$, we will be able to obtain the required estimates and cancelation properties for~$A_\sigma$; see Section~\ref{sec.kernel.properties}.

\medskip

\noindent {\sc Singular and degenerate equations.} The hypotheses made in
Theorem~\ref{th:main} excludes  all the powers
$\varphi(u)=|u|^{m-1}u$ for $m>0$, $m\neq1$, since they are degenerate ($m>1$)
or singular ($m<1$) at the level $u=0$. Nevertheless, a close look at our proof shows
that  we may in fact get  a \lq\lq local'' result,
Theorem~\ref{thm:local.regularity}. Therefore,  we get for these
nonlinearities (and also for more general ones) a    regularity result in the positivity (negativity) set of the solution that implies that bounded weak solutions with a sign are classical; see Section~\ref{sec.degenerate}.

\medskip

\noindent {\sc Higher regularity.} If $\varphi$ is $C^\infty$  we prove that solutions are $C^\infty$. The result will be a consequence of the regularity already provided by Theorem~\ref{th:main} plus a result for linear equations with variable coefficients, Theorem~\ref{th:regularity-linear}, which has an independent interest. The case $\sigma<1$ is a little bit more involved since we first have to raise the regularity in space exponent from $\sigma$ to 1. See more in Section~\ref{sec.higher}, where  higher regularity results depending on the smoothness of $\varphi$ are given.

\begin{Theorem}\label{th:main2}
Let $u$ be a bounded weak solution to equation \eqref{eq:main}. If $\varphi\in C^\infty(\mathbb{R})$, $\varphi'>0$ in $\mathbb{R}$,  then $u\in C^\infty(Q)$.
\end{Theorem}

As a direct precedent of the present
work, let us mention  the paper~\cite{pqrv3}, where we consider the
nonlinearity $\varphi(u)=\log(1+u)$ in the case $\sigma=1=N$, and
prove that solutions with initial data in some $L\log L$ space
become instantaneously bounded and  $C^\infty$. Notice that in this case $\varphi'(u)=1/(1+u)$, and hence the equation is uniformly parabolic.

We expect some of these ideas to have a wider applicability. We point out several possible extensions, together with some comments and applications of equation~\eqref{eq:main} in Section~\ref{sec.extension}.

Let us  remark that Kiselev et al.~give a proof of $C^\infty$ regularity of a class of periodic solutions of geostrophic equations in 2D with $C^\infty$ data~\cite{Kiselev-Nazarov-Volberg}. Their methods are completely different to the ones used in the present paper.

\section{Kernel properties}\label{sec.kernel.properties}
\setcounter{equation}{0}

In this section we consider two issues for the kernel $A_\sigma=(-\Delta)^{\sigma/2}P_\sigma$ which play an important role in the study of regularity, namely some estimates and a cancelation property. Before doing this, it will be convenient to introduce certain H\"older space adapted to equation~\eqref{eq:main}, together with appropriate notations. For simplicity, we will omit the subscript $\sigma$  in what follows when no confusion arises. It will also be convenient to use the notation $Y=(x,t)\in \mathbb{R}^{N+1}$.

The kernel $P$  has as Fourier transform \ $ \widehat
P(\xi,t)=e^{-|\xi|^\sigma t}$. Therefore, it has the self-similar form
\begin{equation}
  \label{kernel}
  P(x,t)=t^{-N/\sigma}\Phi(z),\qquad z=xt^{-1/\sigma}\in\mathbb{R}^N,\quad t>0.
\end{equation}
Moreover, the profile $\Phi$ is a $C^\infty$ positive, radially
decreasing function $\Phi(z)=\widetilde\Phi(|z|)$, satisfying
$\widetilde\Phi(s)\sim s^{-N-\sigma}$ for $s\to\infty$, cf.~\cite{Blumenthal-Getoor}.
We will exploit all these properties fruitfully in what follows.

\medskip

\noindent \textsc{The $\sigma$-distance and the associated H\"older space. }
The self-similar structure of $P$ motivates the use of the {\em $\sigma$-parabolic ``distance''} $|Y_1-Y_2|_\sigma$,  where
\begin{equation*}
\label{sigma-distance}
|Y|_\sigma:=\Big(|x|^2+|t|^{2/\sigma}\Big)^{1/2}=t^{1/\sigma}(|z|^2+1)^{1/2}.
\end{equation*}
This is not really a distance unless $\sigma\ge1$, since the triangle inequality does not hold if $\sigma<1$. However it is a quasimetric, with relaxed triangle inequality
\begin{equation}\label{quasimetric}
|Y_1-Y_3|_\sigma\le 2^{\frac{(1-\sigma)_+}\sigma}\big(|Y_1-Y_2|_\sigma+|Y_2-Y_3|_\sigma\big).
\end{equation}
This will be enough for our purposes.

Observe the relation  between the standard Euclidean distance and this $\sigma$-parabolic distance:
\begin{equation}
\label{eq:relation.metrics}
|Y|\le c|Y|_\sigma^\nu\ \text{ for every }|Y|\le 1,\qquad \nu:=\min\{1,\sigma\}.
\end{equation}

The $\sigma$-parabolic ball  is defined as  $B_R:=\{Y\in\mathbb{R}^{N+1}\,:\,|Y|_\sigma<R\}$.
Performing the change of variables
\begin{equation}\label{ss-change}
s=|x|\,|t|^{-1/\sigma},\quad
r=(|x|^2+|t|^{2/\sigma})^{1/2},
\end{equation}
we get for all $\delta>-N-\sigma$,
$$
\int_{B_R}|Y|_\sigma^{\delta}\,dY=2\sigma N\omega_N\int_0^R
r^{\delta+N+\sigma-1}\,dr\int_0^\infty\frac{
s^{N-1}}{(1+s^2)^{\frac{N+\sigma}2}}\,ds=c R^{\delta+N+\sigma}.
$$
In particular, the volume of the ball  $B_R$  is proportional to
$R^{N+\sigma}$. In the same way, $
\int_{B_R^c}|Y|_\sigma^{-\delta}\,dY=c R^{-\delta+N+\sigma}$ for every $\delta>N+\sigma$.

The  H\"older space  $C^\alpha_\sigma(Q)$, $\alpha\in(0,\nu)$, will consist of functions $u$ defined in $Q$ such that for some constant $c>0$
$$
|u(Y_1)-u(Y_2)|\le c|Y_1-Y_2|_\sigma^\alpha \quad \text{ for every }Y_1,\,Y_2\in Q.
$$

\medskip

\noindent \textsc{The estimates. }
Using formula \eqref{kernel} we deduce  that $A=(-\Delta)^{\sigma/2}P$ has the self-similar expression
\begin{equation*}
  \label{A-selfsimilar}
  A(x,t)= t^{-1-\frac{N}\sigma}\Psi(z),
\end{equation*}
where $z=xt^{-1/\sigma}$, $\Psi(z)=(-\Delta)^{\sigma/2}\Phi(z)$. This is the basis for the estimates.
\begin{Proposition}  \label{prop:derivada de A}
For every $Y\in Q$ the kernel $A$ satisfies
\begin{equation}\label{props-A}
|A(Y)|\le \frac{c}{|Y|_\sigma^{N+\sigma}},\qquad
  |\partial_t A(Y)|\le \frac{c}{|Y|_\sigma^{N+2\sigma}},\qquad|\nabla_x A(Y)|\le \frac{c}{|Y|_\sigma^{N+\sigma+1}}.
\end{equation}
\end{Proposition}

\noindent\emph{Proof. } We observe that
$\widehat\Phi(\xi)=e^{-|\xi|^\sigma}$, hence
$\widehat\Psi(\xi)=|\xi|^{\sigma}e^{-|\xi|^\sigma}$. Using the
expression of the inverse Fourier transform of a radial function, putting $\Psi(z)=\widetilde\Psi(|z|)$,
$$
\widetilde\Psi(s)=c_Ns^{1-\frac N2}\int_0^\infty
e^{-r^\sigma}r^{\frac N2+\sigma} J_{\frac{N-2}2}(rs)\,dr,
$$
we get the decay $|\widetilde\Psi(s)|=O(s^{-N-\sigma})$ for $s$
large \cite[Lemma 1]{Pruitt-Taylor}. Since $\widetilde\Psi$ is bounded, we have
\begin{equation}
\label{eq:estimate.psi}
|\widetilde \Psi(|z|)|=|(-\Delta)^{\sigma/2}\Phi(z)|\le
c(1+|z|^2)^{-\frac{N+\sigma}2},
\end{equation}
which implies
$$
|A(Y)|\le \frac c{t^{1+\frac
N\sigma}(1+|xt^{-1/\sigma}|^2)^{\frac{N+\sigma}2}}=\frac{c}{|Y|_\sigma^{N+\sigma}}.
$$

The estimate for the time derivative is a consequence of
\begin{equation}
\label{eq:relation.A.Phi}
\partial_t A(Y)=-(-\Delta)^\sigma P(x,t)=-t^{-N/\sigma-2}(-\Delta)^{\sigma}\Phi(z),
\end{equation}
which follows immediately from the equation satisfied by $P$,
and~\eqref{eq:estimate.psi}.
Indeed,
\begin{equation*}
\label{eq:estimate.A_t}
|\partial_tA(Y)|\le \frac c{t^{2+\frac
N\sigma}(1+|xt^{-1/\sigma}|^2)^{\frac{N+2\sigma}2}}=\frac{c}{|Y|_\sigma^{N+2\sigma}}.
\end{equation*}

In order to estimate the spatial derivative $\nabla_x A(Y)$, we consider the equation relating the profiles $\Phi$ and $\Psi$,
$$
\sigma(-\Delta)^\sigma\Phi(z)-(N+\sigma)\Psi(z)-z\cdot\nabla\Psi(z)=0,
$$
which follows from~\eqref{eq:relation.A.Phi}. It implies that
$$
|\nabla\Psi(z)|\le \frac c{|z|}(|\Psi(z)|+|(-\Delta)^\sigma\Phi(z)|).
$$
Since $\nabla\Psi$ is bounded, we deduce the estimate $|\widetilde\Psi'(s)|\le c(1+s^2)^{-\frac{N+\sigma+1}2}$. Finally
$$
|\nabla_x A(Y)|= t^{-1-\frac{N+1}\sigma}|\widetilde\Psi'(s)|\le\frac
c{t^{1+\frac{N+1}\sigma}(1+|xt^{-1/\sigma}|^2)^{\frac{N+\sigma+1}2}}=\frac{c}{|Y|_\sigma^{N+\sigma+1}}.
$$
\qed

\medskip

Let us point out that further derivatives may be estimated in a similar way.

\medskip

\noindent \textsc{Cancelation. }
We now show that the function $A$ has zero integral in the sense of principal value  adapted  to the self-similar variables:  we take out a small $\sigma$-ball and integrate, and then  pass to the limit.

\begin{Proposition}\label{pro:average0}
For every $R>\varepsilon>0$,
\begin{equation}
  \label{int-A=0}
\int_{B_R^+-B_\varepsilon}A(x,t)\,dxdt=\int_{B_R^--B_\varepsilon}A(x,t)\,dxdt=0
\end{equation}
where $B_R^+=B_R\cap\{t>0\}$,
$B_R^-=B_R\cap\{t<0\}$.\end{Proposition}

\noindent{\it Proof. } From the equation for the profile $\Phi$,
 $$
 \sigma (-\Delta)^{\sigma/2}\Phi(z)-N\Phi(z)- z\cdot \nabla\Phi(z)=0,
 $$
 we get an alternative expression for the profile of $A$,
\begin{equation*}
\label{eq:alternate-psi} \widetilde\Psi(s)=\frac1\sigma
(N\widetilde\Phi(s)+s\widetilde\Phi'(s))=\frac{
s^{1-N}}{\sigma}(s^N\widetilde\Phi(s))'.
\end{equation*}
Hence, using the change of variables
\eqref{ss-change}  and the behaviour of $\widetilde\Phi$ at infinity, we get
$$
\begin{array}{rl}
\displaystyle\int_{B_R^\pm-B_\varepsilon}A(x,t)\,dxdt&\displaystyle=N\omega_N\sigma
\int_\varepsilon^R\int_0^\infty \frac{s^{N-1}\widetilde\Psi(s)}{r}\,ds dr \\ [3mm]
&\displaystyle=\left.N\omega_N\log(
R/\varepsilon)(s^N\widetilde\Phi(s))\right|_0^\infty=0.
\end{array}
$$
 \qed

\section{The linear problem}
\label{sec.linear}
\setcounter{equation}{0}
As we have said in the Introduction, the solution $u$ to equation~\eqref{eq:main} will be analyzed by writing it as a solution of a linear problem with a particular right hand side. This leads to the representation of $u$ by  means of a variation of constants formula. We give a  proof of this independent fact  and then proceed to establish the regularity of this linear problem.

\subsection{A representation formula}\label{sec.rep}

Let us consider the Cauchy problem associated to the fractional linear
heat equation with a source term,
\begin{equation}  \label{eq:linear-con-f}
\left\{
\begin{array}{ll}
\partial_t u+(-\Delta)^{\sigma/2}u=(-\Delta)^{\sigma/2}f,\qquad & (x,t)\in  Q,
\\ [4mm]
u(x,0) = u_0(x),\qquad & x\in\mathbb{R}^N.%
\end{array}
\right.
\end{equation}
We assume that $f\in L_{\rm loc}^\infty((0,\infty):L^1(\R^N,\rho\,dx))$ with  $\rho(x)=(1+|x|)^{-(N+\sigma)}$. We define a  {\sl very weak solution} to problem \eqref{eq:linear-con-f}  as a function   $u\in C([0,\infty): L^1(\R^N,
\rho\,dx))$, such that
\[
\iint_{Q} u(x,t)\partial_t\zeta(x,t)\,dx\,dt =\int_{Q} (u-f)(x,t)(-\Delta)^{\sigma/2}\zeta(x,t)\,dx\,dt
\]
for all $\zeta\in C_c^\infty(Q)$, and $u(x,0)=u_0$ almost everywhere.

\begin{Theorem}
\label{th:duhamel.very.weak}   If $f\in L_{\rm loc}^\infty([0,\infty):L^1(\R^N,\rho\,dx))\cap
C^\alpha_\sigma(Q)$,  $0<\alpha\le\min\{1,\sigma\}$, and $u_0\in
L^p(\mathbb{R}^N)$ for some $1\le p\le\infty$, there is a unique
very weak solution of problem~\eqref{eq:linear-con-f}, which is given by
representation using Duhamel's formula:
\begin{equation}\label{eq:duhamel1}
\begin{array}{l}
u(x,t)=\displaystyle\int_{\mathbb{R}^N}P(x-\overline x, t) u_0(\overline x)\,d\overline x\\[4pt]
\qquad\qquad\qquad \displaystyle+\int_0^t\int_{\mathbb{R}^N}(-\Delta)^{\sigma/2}P(x-\overline x,
t-\overline t) f(\overline x,\overline t)\,d\overline xd\overline t,\quad (x,t)\in Q.
\end{array}
\end{equation}
\end{Theorem}

\noindent{\it Proof. } \noindent\textsc{Step 1. Uniqueness. } We may
assume $u_0=f=0$ and then apply the results of \cite{Barrios-Peral-Soria-Valdinoci} where a wider class of data and solutions is treated.

\noindent\textsc{Step 2. $u$ is well defined. } We only have to take care of the last term in~\eqref{eq:duhamel1}, which can be written as
\begin{equation*}
\label{eq:duhamel2}
\int_QA(Y-\overline Y)\chi_{\{\overline t<t\}} f(\overline Y)\,d\overline Y,
\end{equation*}
$Y\in Q\subset \mathbb{R}^{N+1}$.
In order to prove that this integral is well defined we decompose $Q=B^-_r\cup (Q\backslash B^-_r)$, with $r>0$ small, where
$B^-_r=\{\overline Y=(\overline x,\overline t)\,:\,|Y-\overline Y|_\sigma<r, \ \overline t\le t\}$. The cancellation property \eqref{int-A=0} allows us to estimate the inner integral,
$$
\begin{array}{rcl}
\displaystyle\left|\int_{B^-_r}A(Y-\overline Y)\chi_{\{\overline t<t\}}f(\overline Y)\,d\overline Y\right|&=&
\displaystyle
\left|\int_{B^-_r}A(Y-\overline Y)\chi_{\{\overline t<t\}}\Big(f(\overline Y)-f(Y)\Big)\,d\overline Y\right|
\\ [4mm]
&\le&\displaystyle\int_{B^-_r}|A(Y-\overline Y)|\,|f(\overline Y)-f(Y)|\,d\overline Y \\
[4mm] &\le&\displaystyle
c\int_{B^-_r}\frac{d\overline Y}{|\overline Y- Y|_\sigma^{N+\sigma-\alpha}}\le c.
\end{array}
$$
The outer integral is bounded by using estimate \eqref{props-A}.

\noindent\textsc{Step 3. $u$ is a very weak solution. }
In order to justify the representation formula we proceed by approximation. Let $t>0$  and take $f\in C^\infty_c(Q)$
with $f(x,\overline t)=0$ for $\overline t\ge t-r$, $r$ small, thus avoiding the singularity.  In that case the integral in the ball $B^-_r$ vanishes identically. Since
moreover $u$ given by \eqref{eq:duhamel1} is a bounded classical solution, hence a very weak solution, the assertion holds.

Next, for any  $f\in C^\infty(Q)$ and compactly supported in space (in a uniform way), we use approximation with functions $f_n$ as before by modifying $f$ in the time interval
$t-r_n\le \overline t\le t$. Using the fact that the fractional Laplacian can be applied to $f_n$ instead of $P$, it is easy to see that we can pass to the limit $u$ of the solutions $u_n$, which is still a bounded
classical solution. Moreover, the formula as it is written holds for functions
$f$ of this class by integrating by parts and the integrability estimate from Step 2.

Finally, for general $f$ as in the statement, we use approximation of $f$ in
a compact set by functions $f_n\in C^\infty(Q)$ compactly supported in space. Passing to the
limit in the very weak formulation, which is again justified thanks to Step 2, we obtain that $u=\lim u_n$ is a very weak
solution.  \qed

\subsection{Regularity of the linear problem}\label{sec.reg.linear}

The first term in the right-hand side
of~the representation formula~\eqref{eq:duhamel1} is regular. Hence, $u$ has the same regularity as
\begin{equation}
  \label{eq:g}
  g(Y)=\int_{\mathbb{R}_+^{N+1}}A(Y-\overline Y)\chi_{\{\overline t<t\}}f(\overline Y)\,d\overline Y.
\end{equation}
We start by proving that $g$ has the same $\sigma$-H\"older regularity as $f$.
\begin{Lemma}
\label{lem:regularity.improved}
Let $f\in C^{\alpha}_\sigma(Q)\cap
L^\infty(Q)$ for some $0<\alpha<\nu$, and $g$ given by~\eqref{eq:g}. Then
 $g\in
C^{\alpha}_\sigma(Q)\cap L^\infty(Q)$.
\end{Lemma}

\begin{figure}[ht]
\begin{center}
\epsfig{file=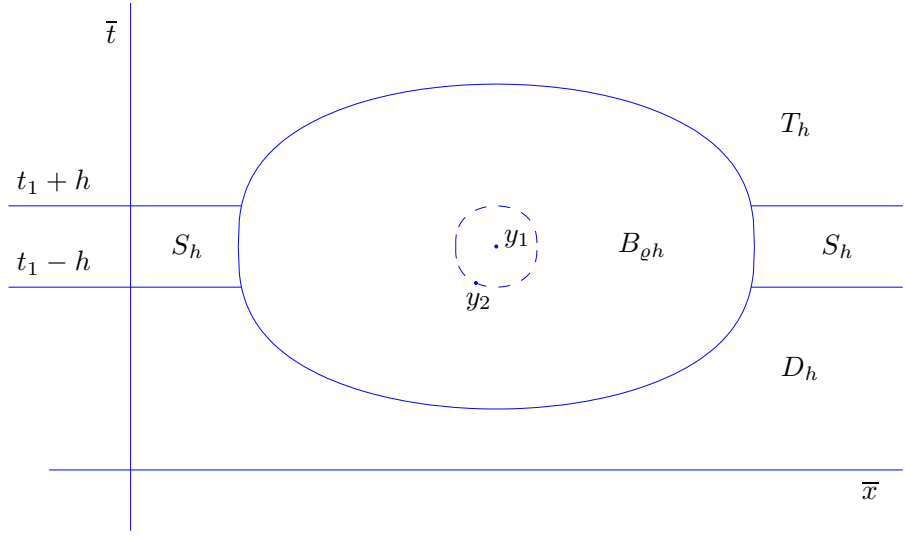,width=12cm}
\end{center}
\caption{Integration regions.}
\label{fig:integration.regions}
\end{figure}

\noindent{\it Proof. }  Let $Y_1=(x_1,t_1),\,Y_2=(x_2,t_2)\in Q$ be two points with $|Y_1-
Y_2|_\sigma=h>0$ small. We have to estimate the difference
\begin{equation}\label{eq:g-g0}
g(Y_1)-g(Y_2)=\int_{\mathbb{R}_+^{N+1}}\Big(A( Y_1-
\overline Y)\chi_{\{\overline t<t_1\}}-A( Y_2- \overline Y)\chi_{\{\overline t<t_2\}}\Big)f(\overline Y)\,d\overline Y
\end{equation}
and see if we get  that it is $O(h^\alpha)$. We decompose $Q$ into
four regions, depending on the sizes of $|\overline x-x_1|$ and $|\overline t-t_1|$, see
Figure~\ref{fig:integration.regions}.

\medskip

\noindent (i) \emph{The small ``ball'' $B_{\varrho h}(Y_1)$}, where $\varrho>1$ is a constant to be fixed later.  We take $h$ small enough ($\varrho h<\min\{t_1,1\}$) so that, on the one hand
$B_{\varrho h}\subset Q$, and, on the other hand, we can use the relation~\eqref{eq:relation.metrics}.   The difficulty in this region is the
non-integrable singularity of $A( Y)$ at $ Y=0$. Integrability will
be gained thanks to the regularity of $f$.  We first have, repeating the computations in Step~2 of the proof of Theorem~\ref{th:duhamel.very.weak},
$$
\displaystyle\left|\int_{B_{\varrho h(Y_1)}}A( Y_1-
\overline Y)\chi_{\{\overline t<t_1\}}f(\overline Y)\,d\overline Y\right|\displaystyle
\le
ch^{\alpha}.
$$
To estimate the second term in \eqref{eq:g-g0}, we consider the ball
$B_h(Y_2)$. To be sure that the distance from $\partial B_{\varrho h}(Y_1)$ to $B_h(Y_2)$ is positive, we take $\varrho=\max\{4,2^{2/\sigma}\}$; see \eqref{quasimetric}. Using again
the cancelation property \eqref{int-A=0}, we get
\begin{equation}\label{integral-I1}
\begin{array}{l}
\displaystyle\int_{B_{\varrho h}(Y_1)}A(Y_2-\overline Y)\chi_{\{\overline t<t_1\}}f(\overline Y)\,d\overline Y=
\\[4mm]
\qquad\qquad\underbrace{\displaystyle\int_{B_{\varrho h}(Y_1)}A(Y_2-\overline Y)\chi_{\{\overline t<t_2\}}
\Big(f(\overline Y)-f(Y_2)\Big)\,d\overline Y}_{I_1}\\[4mm]
\qquad\qquad\qquad+
\underbrace{f(Y_2)\int_{B_{\varrho h}(Y_1)-B_{h}(Y_2)}A(Y_2-\overline Y)\chi_{\{\overline t<t_2\}}\,d\overline Y}_{I_2}.
\end{array}
\end{equation}
The first integral $I_1$  satisfies, as before, $|I_1|\le ch^\alpha$.
As to $I_2$, since we are far from the singularity of $A$,
$$
|I_2|\le c
h^{\alpha}\int_{B_{\varrho h}(Y_1)-B_{h}(Y_2)}\frac{d\overline Y}{h^{N+\sigma}}\le
ch^{\alpha}.
$$

\noindent (ii) \emph{The narrow strip
$S_h=\{\overline Y\in \overline B_{\varrho h}^{\; c}(Y_1),\;|\overline t-t_1|<h^\sigma\}$. } In this region
we have $|\overline Y -Y_1|_\sigma\le \varrho_1|\overline Y -Y_2|_\sigma$ and
$|\overline x-x_1|>\varrho_2h$, for some positive constants $\varrho_1,\,\varrho_2$ depending only on $\sigma$. Therefore,
$$
\begin{array}{l}
\left|\displaystyle\int_{S_h}
\Big(A(Y_1-\overline Y)\chi_{\{\overline t<t_1\}}-A(Y_2-\overline Y)\chi_{\{\overline t<t_2\}}\Big)f(\overline Y)\,d\overline Y\right|\\[4mm]
\qquad\displaystyle\le
\int_{S_h}\Big(|A(Y_1-\overline Y)|+|A(Y_2-\overline Y)|\Big)|f(\overline Y)|\,d\overline Y \le
\int_{S_h}\frac{d\overline Y}{|\overline Y -Y_1|_\sigma^{N+\sigma-\alpha}}\\[4mm]
\qquad\displaystyle\le
c\int_{t_1-h^\sigma}^{t_1+h^\sigma}\int_{\{|\overline x-x_1|>\varrho_2h\}}\frac{d\overline xd\overline t}{|\overline x-x_1|^{N+\sigma-\alpha}}
\displaystyle\le ch^{\alpha}.
\end{array}
$$
Notice that   $\alpha<\sigma$, so that the last integral is convergent.

\noindent (iii) \emph{The complement of the ball $B_{\varrho h}(Y_1)$ for large
times, $T_h=\{\overline Y\in \overline B_{\varrho h}^{\; c}(Y_1),\;\overline t>t_1+h^\sigma\}$. } The
integral in this region is 0, since here we have
$$
A(Y_1-\overline Y)\chi_{\{\overline t<t_1\}}=A(Y_2-\overline Y)\chi_{\{\overline t<t_2\}}=0.
$$

\medskip

\noindent (iv) \emph{The complement of the ball $B_{\varrho h}(Y_1)$ for small
times, $D_h=\{\overline Y\in \overline B_{\varrho h}^{\; c}(Y_1),\;\overline t<t_1-h^\sigma\}$. } The
required estimate is obtained here using the fact that
$$
A(Y_1-\overline Y)\chi_{\{\overline t<t_1\}}-A(Y_2-\overline Y)\chi_{\{\overline t<t_2\}}=A(Y_1-\overline Y)-A(Y_2-\overline Y).
$$
Thus we are
integrating  a difference of $A$'s, so there will be some
cancelation. Indeed, by the Mean Value Theorem,
$$
|A(Y_1-\overline Y)-A(Y_2-\overline Y)|\le |Y_1-Y_2| \max\{|\nabla_x A(\theta)|, |\partial_t
A(\theta)|\},
$$
where $\theta=s(Y_1-\overline Y)+(1-s)(Y_2-\overline Y)$ for some $s\in(0,1)$.

Since $|Y_1-Y_2|\le|Y_1-Y_2|_\sigma^\nu=h^\nu$, and in $D_h$ it holds $|Y_1-\overline Y|_\sigma\le2|\theta|_\sigma$, Proposition~\ref{prop:derivada de A}
yields
$$
|A(Y_1-\overline Y)-A(Y_2-\overline Y)|\le
\frac{ch^\nu}{|\theta|_\sigma^{N+\sigma+\nu}}\le\frac{ch^\nu}{|\overline Y- Y_1|_\sigma^{N+\sigma+\nu}}.
$$
Therefore we conclude that
$$
\begin{array}{l}
\displaystyle\left|\int_{D_h}\Big(A(Y_1-\overline Y)\chi_{\{\overline t<t_1\}}-A(Y_2-\overline Y)\chi_{\{\overline t<t_2\}}\Big)f(\overline Y)\,d\overline Y\right| \\
\hspace{25mm}\displaystyle\le
ch^\nu\int_{D_h}\frac{d\overline Y}{|Y_1-\overline Y|_\sigma^{N+\sigma+\nu-\alpha}}\le
ch^\nu\int_{\{|\overline Y|_\sigma>\varrho h\}}\frac{d\overline Y}{|\overline Y|_\sigma^{N+\sigma+\nu-\alpha}}=
ch^{\alpha},
\end{array}
$$
assuming  $\alpha<\nu$ so that the last integral is convergent.
\qed

\noindent\emph{Remark. } Notice that if some derivative (even a fractional one) of $f$ belongs to $C^\alpha_\sigma(Q)\cap L^\infty(Q)$, then a computation analogous to that in the above lemma shows that the convolution of this derivative against the kernel $A$ also belongs to $C^\alpha_\sigma(Q)\cap L^\infty(Q)$. We conclude that $g$ has the same regularity as $f$.

\medskip

As a corollary of Lemma~\ref{lem:regularity.improved}, we obtain a maximal regularity result for the linear equation with a standard right hand side that has independent interest.

\begin{Corollary} Let $f\in C^{\alpha}_{\sigma}(Q)\cap
L^\infty(Q)$, $0<\alpha<\min\{1,\sigma\}$. If $w$ is a very weak solution to
\begin{equation}
\label{eq:w} \partial_tw+(-\Delta)^{\sigma/2} w=f,
\end{equation}
then $\partial_t w, (-\Delta)^{\sigma/2} w\in C^{\alpha}_{\sigma}(Q)\cap
L^\infty(Q)$.
\end{Corollary}

\noindent\emph{Proof. } The function $u=(-\Delta)^{\sigma/2} w$ solves~\eqref{eq:linear-con-f} in the distributional sense. Hence $u\in C^\alpha_\sigma(Q)\cap
L^\infty(Q)$. The result now follows noticing that $\partial_t w=f-u$.
\qed

\section{Improving $\sigma$-H\"older regularity}
\label{sec.regularity.nonlinear}
\setcounter{equation}{0}

We now return to the nonlinear equation~\eqref{eq:main}. For bounded weak energy solutions the equation is neither degenerate nor
singular. Hence, the results from~\cite{Athanasopoulos-Caffarelli} guarantee that they are $C^\alpha_\sigma$ for some
$\alpha\in(0,\nu)$. The aim of this section is to improve this regularity showing that the solutions belong to $C_\sigma^\alpha$ for \emph{all} $\alpha<\nu$. Further regularity, showing that the solution is classical, will be obtained in Section~\ref{sec.beyond.Holder}.

The idea is to use that the solution $u$ to the nonlinear equation~\eqref{eq:main} is a solution to the linear equation
$$
\partial_t u+(-\Delta)^{\sigma/2}u=(-\Delta)^{\sigma/2}f,\qquad  f(Y)=u(Y)-\frac{\varphi(u(Y))}{\varphi'(u(Y_0))}.
$$
Since $u\in C_\sigma^\alpha(Q)$, $\varphi$ is uniformly parabolic,  and $\varphi'\in C^{\gamma}(\mathbb{R})$, applying the Mean Value Theorem we get
\begin{equation}
\label{eq:local.holder.condition0}
\begin{array}{rcl}
|f(Y_1)-f(Y_2)|&=&\displaystyle
\left|1-\frac{\varphi'(\theta)}{\varphi'(u(Y_0))}\right|\,|u(Y_1)-u(Y_2)|\\[4mm]
&\le&\displaystyle\frac{|u(Y_0)-\theta|^\gamma}{\varphi'(u(Y_0))}\,|Y_1-Y_2|_\sigma^{\alpha} \\ [4mm]
&\le&\displaystyle
c\max\{|u(Y_1)-u(Y_0)|^\gamma,|u(Y_2)-u(Y_0)|^\gamma\}|Y_1-Y_2|_\sigma^{\alpha}\\
[4mm] &\le&\displaystyle
c\max\{|Y_1-Y_0|_\sigma^{\alpha\gamma},|Y_2-Y_0|_\sigma^{\alpha\gamma}\}|Y_1-Y_2|_\sigma^{\alpha},
\end{array}
\end{equation}
where $\theta$ is some value between $u(Y_1)$ and $u(Y_2)$. This gives not only that $f$ has the same regularity as $u$, namely  $f\in C^{\alpha}_\sigma(Q)$, but a bit more that will be enough to improve the $\sigma$-H\"older regularity of $u$ by a constant factor.

\begin{Lemma}
\label{lem:regularity.improved2} Let $f\in L^\infty(Q)$ and let $g$ be the function defined in \eqref{eq:g}. Assume
that there exist  $c>0$, $\delta_0>0$ and $\epsilon>0$,
$\alpha+\epsilon<\nu$,  such that
\begin{eqnarray}
\label{eq:local.holder.condition2}
|f(Y)-f(Y_0)|\le c|Y-Y_0|_\sigma^{\alpha+\epsilon},\\ [3mm]
\label{eq:local.holder.condition}
|f(Y_1)-f(Y_2)|\le c\delta^{\epsilon}\,|Y_1-Y_2|_\sigma^{\alpha},
\end{eqnarray}
for all $0<\delta<\delta_0$, $Y,Y_1,Y_2\in
B_\delta(Y_0)$. Then,
\begin{equation*}
\label{eq:g.alpha+eps} |g(Y)-g(Y_0)|\le
c'|Y-Y_0|_\sigma^{\alpha+\epsilon},
\end{equation*}
for all $Y\in
B_\delta(Y_0)$, where $c'$ depends on $c$.
\end{Lemma}

\noindent{\it Proof. }  The fact that $f$ is
$C^{\alpha+\epsilon}_\sigma$ at $Y_0$ (with $\alpha+\epsilon<\nu$) implies that all the estimates used to prove Lemma~\ref{lem:regularity.improved}  work and yield terms which
are $O(h^{\alpha+\epsilon})$, except that for the integral $I_1$ in
\eqref{integral-I1}. To estimate this term take $\varrho h<\delta_0$ and observe that~\eqref{eq:local.holder.condition} gives
$$
\begin{array}{rl}
|I_1|&=\displaystyle\left|\int_{B_{\varrho h}(Y_0)}A(Y-\overline Y)\chi_{\{\overline t<t\}}
\Big(f(\overline Y)-f(Y)\Big)\,d\overline Y\right| \\ [4mm]
&\displaystyle\le
ch^{\epsilon}\int_{B_{\varrho h}(Y_0)}\frac1{|Y-\overline Y|_\sigma^{N+\sigma}}|Y-\overline Y|_\sigma^\alpha\,d\overline Y\le
ch^{\alpha+\epsilon}.
\end{array}
$$
\qed

Applying this lemma a finite number of times we obtain the desired regularity.

\begin{Theorem}\label{cor:todo.alpha}  Let $u$ be a bounded weak solution to equation \eqref{eq:main}. Assume $\varphi\in C^{1,\gamma}(\mathbb{R})$  and $\varphi'(s)>0$ for every $s\in\mathbb{R}$. Then $u$ belongs to
$C^\alpha_\sigma(Q)$ for all $\alpha\in(0,\nu)$.
\end{Theorem}

We must remark that the restriction $\alpha+\epsilon<\nu$ in Lemma~\ref{lem:regularity.improved2} is only needed to make the outer integrals convergent; the estimate in the ball $B_{\varrho h}(Y_0)$ is true for any $\alpha\in(0,\nu)$, $\epsilon\in(0,1]$. This observation turns to be of great importance in obtaining further regularity in the next section.

\section{Classical solutions}
\label{sec.beyond.Holder}
\setcounter{equation}{0}

Our  next aim is to go beyond the $C^\nu_\sigma$ threshold of regularity. We encounter here an additional difficulty, steaming from the nonlocal  character of the fractional Laplacian operator,  which is not present in the work~\cite{Caffarelli-Vasseur}, namely that $A(Y-\overline Y)\chi_{\{\overline t<t\}}\neq(-\Delta)^{\sigma/2}(P(Y-\overline Y)\chi_{\{\overline t<t\}})$. For that reason we must treat the second order estimates in the time and space variables separately. We begin by improving regularity in space, to
obtain $u(\cdot,t)\in C^{\alpha}(\mathbb{R}^N)$ uniformly in $t$ for some $\alpha>\nu$ depending on the regularity of the nonlinearity $\varphi$. We then use equation~\eqref{eq:main} to get Lipschitz regularity in time, which is later improved to get  $u(x,\cdot)\in C^{\nu(1+\gamma)/\sigma}(\mathbb{R}_+)$ uniformly in $x$. The last step is to reach  the desired smoothness in space, $u(\cdot,t)\in C^{\nu(1+\gamma)}(\mathbb{R}^N)$ uniformly in $t$.

\noindent\emph{Notation. } By $u\in C^\alpha$ with  $\alpha\in [1,2)$ we mean $u\in C^{1,\alpha-1}$ if $\alpha\in(1,2)$, and $u\in C^{0,1}$ if $\alpha=1$.

\begin{Lemma}
\label{lem:regularity.improved3} Let $f\in L^\infty(Q)$ satisfy \eqref{eq:local.holder.condition2} and \eqref{eq:local.holder.condition} with $0<\alpha<\nu$, $0<\epsilon<1$, and let $g$ be the function defined in \eqref{eq:g}. Then, for every $e\in\mathbb{R}^N$, $|e|=1$,
\begin{equation}
\label{eq:g.alpha+eps+mas}
|g(x_0+he,t_0)-2g(x_0,t_0)+g(x_0-he,t_0)|\le ch^{\alpha+\epsilon}\quad \text{for every $h>0$ small}.
\end{equation}
\end{Lemma}

\noindent\emph{Proof. }
Put $Y=Y_0+(he,0)$ and let $Y^*=2Y_0-Y$ be its symmetric point with respect to $Y_0$. We have to estimate the second difference
$$
g(Y)-2g(Y_0)+g(Y^*) =\int_{\mathbb{R}_+^{N+1}}\mathcal{A}(Y,Y_0,\overline Y)f(\overline Y)\,d\overline Y,
$$
where
$$
\begin{array}{rl}
\displaystyle\mathcal{A}(Y,Y_0,\overline Y)&\displaystyle=A(Y-\overline Y)\chi_{\{\overline t<t\}}-2A(Y_0-\overline Y)
\chi_{\{\overline t<t_0\}}+A(Y^*-\overline Y)\chi_{\{\overline t<t^*\}} \\ [3mm]
&\displaystyle=\Big(A(Y-\overline Y)-2A(Y_0-\overline Y)
+A(Y^*-\overline Y)\Big)\chi_{\{\overline t<t_0\}}.
\end{array}
$$
As in the  proof of Lemma~\ref{lem:regularity.improved}, we consider
separately the contributions to the integral in several regions, though here we only need to consider the ball $B_{\varrho h}(Y_0)$ and its complement $B_{\varrho h}^c(Y_0)$.
The contribution of the integral in $B_{\varrho h}(Y_0)$ is decomposed
as the sum $J_1-2J_2+J_3$, where
$$
\begin{array}{l}
\displaystyle J_1=\int_{B_{\varrho h}(Y_0)} A(Y-\overline Y)\chi_{\{\overline t<t_0\}}f(\overline Y)\,d\overline Y,\\[10pt]
\displaystyle J_2=\int_{B_{\varrho h}(Y_0)} A(Y_0-\overline Y)\chi_{\{\overline t<t_0\}}f(\overline Y)\,d\overline Y,\\[10pt]
\displaystyle J_3=\int_{B_{\varrho h}(Y_0)} A(
Y^*-\overline Y)\chi_{\{\overline t<t_0\}}f(\overline Y)\,d\overline Y.
\end{array}
$$
Thus $|\int_{B_{\varrho h}(Y_0)}\mathcal{A}f|\le |J_1|+2|J_2|+|J_3|$. The integrals $J_1$ and $J_2$ were already estimated in the course
of the  proof of Lemma~\ref{lem:regularity.improved}, modified as in Lemma~\ref{lem:regularity.improved2}  (see the comment after that lemma), thus obtaining $O(h^{\alpha+\epsilon}$). Since
 $Y^*-Y_0= Y_0-Y$, the integral $J_3$  is
estimated just in the same way.

To estimate the contribution in the complement of the ball we use Taylor's formula. We have, by Proposition~\ref{prop:derivada de A},
$$
\begin{array}{rl}
\displaystyle|\mathcal{A}(Y,Y_0,\overline Y)|&\displaystyle= |A(Y-\overline Y)-2A(Y_0-\overline Y)+A(2Y_0-Y-\overline Y)| \\ [3mm]
&\displaystyle\le ch^2|D^2_x A(\theta)|\le
\frac{ch^2}{| Y_0-\overline Y|_\sigma^{N+\sigma+2}},
\end{array}
$$
where $\theta$ is as before some intermediate point. This gives
$$
\int_{B_{\varrho h}^c(Y_0)} |\mathcal{A}(Y,Y_0,\overline Y)|\,|f(\overline Y)|\,d\overline Y \le c h^2
\int_{B_{\varrho h}^c(Y_0)}\frac{d\overline Y}{|Y_0-\overline Y|_\sigma^{N+\sigma+2-\alpha-\epsilon}} =
ch^{\alpha+\epsilon}.
$$
We have used that $\alpha+\epsilon<2$, and so the integral converges. This completes the desired estimate. \qed

\begin{Lemma}\label{th:submain1} Under the hypotheses of Theorem~{\rm \ref{th:main}},  bounded weak solutions $u$ to equation \eqref{eq:main} satisfy
$u(\cdot,t)\in C^{\alpha(1+\gamma)}(\mathbb{R}^N)$ for every $\alpha\in(0,\nu)$
uniformly in $t\ge\tau>0$ for every $\tau>0$.
\end{Lemma}

\noindent\emph{Proof. } For each given $Y_0\in Q$ we define a function $g$ and,  as in the proof of Lemma~\ref{lem:regularity.improved2}, we deduce estimate~\eqref{eq:g.alpha+eps+mas} with $\epsilon=\alpha\gamma$ at that point, which is translated into the same estimate for the solution $u$. Since the constants do not depend on the particular point chosen, we get that $u$ satisfies
$$
|u(x+he,t)-2u(x,t)+u(x-he,t)|\le ch^{\alpha(1+\gamma)},
$$
with constant uniform in $Q$. We can thus prove  that
$(-\Delta)^{\delta/2} u$ is bounded in $\mathbb{R}^N$ for every $t>\tau>0$ and  every
$\delta\in(0,\alpha(1+\gamma))$. Indeed,
\begin{equation}\label{lapa-bdd}
\begin{array}{rl}
\displaystyle|(-\Delta)^{\delta/2}u(x,t)|&\displaystyle=
\left|c_\delta\int_{\mathbb{R}^{N+1}}\frac{u(x+z,t)-2u(x,t)+u(x-z,t)}{|z|^{N+\delta}}dz\right|
\\ [4mm]&\displaystyle\le
c\int_{\{|z|<\tau\}}\frac{|z|^{{\alpha(1+\gamma)}}}{|z|^{N+\delta}}dz+
c\int_{\{|z|>\tau\}}\frac{dz}{|z|^{N+\delta}} \le c.
\end{array}
\end{equation}
The result now follows from~\cite[Proposition
2.9]{Silvestre-2007}.
\qed

\begin{Lemma}\label{th:submain2} Under the hypotheses of Theorem~{\rm \ref{th:main}},  $u(x,\cdot)\in
C^{\nu(1+\gamma)/\sigma}(\mathbb{R}_+)$ uniformly in $x\in\mathbb{R}^N$.
\end{Lemma}

\noindent\emph{Proof. }
We first show that $|(-\Delta)^{\sigma/2}\varphi(u)|$ is bounded in $Q$. For that purpose we estimate the second differences in $x$ of $\varphi(u)$ in terms of second differences in $x$ of $u$ and use the previous result. If $Z=(he,0)$, $e\in\mathbb{R}^N$, $|e|=1$,
$$
\begin{array}{l}
\displaystyle|\varphi(u(Y_0+Z))-2\varphi(u(Y_0))+\varphi(u(Y_0-Z))| \\ [3mm]
\displaystyle\qquad\qquad\le|\varphi(u(Y_0+Z))-2\varphi(u(Y_0))+\varphi(2u(Y_0)-u(Y_0+Z))|\\ [3mm]
\displaystyle\qquad\qquad\qquad+|\varphi(u(Y_0-Z))-\varphi(2u(Y_0)-u(Y_0+Z))|\\ [3mm]
\displaystyle\qquad\qquad\le [\varphi]_{C^{1+\gamma}}|u(Y_0+Z)-u(Y_0)|^{1+\gamma}\\ [3mm]
\displaystyle\qquad\qquad\qquad+\|\varphi'(u)\|_{L^\infty(Q)}|u(Y_0+Z))-2u(Y_0)+u(Y_0-Z)|\le ch^{\alpha(1+\gamma)}
\end{array}
$$
for every $\alpha<\nu$. Since $\nu(1+\gamma)>\sigma$ we get, analogously to how we obtained \eqref{lapa-bdd}, that $|(-\Delta)^{\sigma/2}\varphi(u)|\le c$ in $Q$.
Now, using  the equation we get  that $|\partial_tu|\le c$ in $Q$, that is, $u$ is Lipschitz continuous in time, uniformly in space. This means $u\in C^\sigma_\sigma(Q)$.
With this information we now try to repeat the above calculations of Lemma~\ref{th:submain1} with $x_0$ fixed and varying $t$. To this end we consider the point $Y=Y_0+(0,h)$, $h>0$ (for simplicity), and we replace $h$ by $h^{1/\sigma}$ in the regions of integration, see the proof of Theorem~\ref{lem:regularity.improved}.

First, the integral in the ball $B_{\varrho h^{1/\sigma}}(Y_0)$ is estimated as in Lemma~\ref{lem:regularity.improved}, taking note of \eqref{eq:local.holder.condition0}, which  holds with $\alpha=\nu$. Thus
$|\int_{B_{\varrho h^{1/\sigma}}(Y_0)}{\cal A}f|=O(h^{\nu(1+\gamma)/\sigma}).$

Now consider the region $D_{h^{1/\sigma}}=\{\overline Y\in \overline B_{\varrho h^{1/\sigma}}^{\; c}(Y_0),\;\overline t<t_0-h\}$. The idea here is that the characteristic functions take all the value one. Thus, by using Taylor's expansion, since only $t$ varies, we have
$$
\begin{array}{rl}
\displaystyle|\mathcal{A}(Y,Y_0,\overline Y)|&\displaystyle= |A(Y-\overline Y)-2A(Y_0-\overline Y)+A(2Y_0-Y-\overline Y)| \\ [3mm]
&\displaystyle\le ch^{2}|\partial^2_t A(\theta)|\le
\frac{ch^{2}}{| Y_0-\overline Y|_\sigma^{N+3\sigma}},
\end{array}
$$
where $\theta$ is some intermediate point.
This gives
$$
\int_{ D_{h^{1/\sigma}}} |\mathcal{A}(Y,Y_0,\overline Y)|\,|f(\overline Y)|\,d\overline Y \le c h^{2}
\int_{ D_{h^{1/\sigma}}}\frac{|Y_0-\overline Y|_\sigma^{\nu(1+\gamma)}}{|Y_0-\overline Y|_\sigma^{N+3\sigma}}d\overline Y
\le ch^{\nu(1+\gamma)/\sigma}.
$$

We now turn our attention to the difficult part, the small slice $ S_{h^{1/\sigma}}=\{\overline Y\in \overline B_{\varrho h^{1/\sigma}}^{\; c}(Y_0),\;|\overline t-t_0|<h\}$, where we have to look more carefully at the possible cancelations. We have
$$
\begin{array}{l}
\displaystyle\int_{ S_{h^{1/\sigma}}} \mathcal{A}(Y,Y_0,\overline Y)f(\overline Y)\,d\overline Y
\\[6mm]
\qquad\displaystyle=\int_{ S_{h^{1/\sigma}}} \big(A(Y-\overline Y)\chi_{\{\overline t<t\}}-2A(Y_0-\overline Y)\chi_{\{\overline t<t_0\}}+A(Y^*-\overline Y)\chi_{\{\overline t<t^*\}}\big)f(\overline Y)\,d\overline Y
 \\[6mm]
\qquad\displaystyle=\int_{ S_{h^{1/\sigma}}} \big(A(Y-\overline Y)-2A(Y_0-\overline Y)\chi_{\{\overline t<t_0\}}\big)f(\overline Y)\,d\overline Y
\\[6mm]
\qquad\displaystyle=
 \underbrace{\int_{ S_{h^{1/\sigma}}} \big(A(Y-\overline Y)-A(Y_0-\overline Y)\big)f(\overline Y)\,d\overline Y}_{J_1}
\\[6mm]
\qquad\displaystyle\qquad\qquad+\underbrace{\int_{ S_{h^{1/\sigma}}} \big(A(Y_0-\overline Y)\chi_{\{\overline t>t_0\}}-A(Y_0-\overline Y)\chi_{\{\overline t<t_0\}}\big)f(\overline Y)\,d\overline Y}_{J_2}.
\end{array}
$$

First, by the Mean Value Theorem applied to $A$ in the time variable, together with the regularity $C_\sigma^\nu$ of $u$ and Lemma~\ref{lem:regularity.improved}, we have
$$
\begin{array}{rl}
|J_1|&\displaystyle\le\int_{t_0-h}^{t_0+h}\int_{\{|\overline x-x_0|>\varrho h^{1/\sigma}\}}
\frac{ch|\overline Y-Y_0|_\sigma^{\nu(1+\gamma)}}{|\overline Y-Y_0|_\sigma^{N+2\sigma}}\,d\overline xd\overline t \\ [6mm]
&\displaystyle\le ch^2\int_{\{|\overline x-x_0|>\varrho h^{1/\sigma}\}}
\frac{d\overline x}{|\overline x-x_0|^{N+2\sigma-\nu(1+\gamma)}}
=ch^{\nu(1+\gamma)/\sigma}.
\end{array}
$$
As to the second integral $J_2$,  performing the change of variables $\overline Y\to Z_1= \overline Y^*\equiv (\overline x,2t_0-\overline t)$, symmetric in time, in the second term (and writing again $\overline Y$ instead of $Z_1$), we have
$$
\begin{array}{rl}
J_2&\displaystyle=\int_{S_{h^{1/\sigma}}} A(Y_0-\overline Y)\chi_{\{\overline t>t_0\}} f(\overline Y)\,d\overline Y -\int_{S_{h^{1/\sigma}}} A(\overline Y-Y_0)\chi_{\{\overline t>t_0\}}f(\overline Y^*)\,d\overline Y \\ [6mm]
&\displaystyle=\int_{ S_{h^{1/\sigma}}} A(Y_0-\overline Y)\chi_{\{\overline t>t_0\}}\big(f(\overline Y)-f(\overline Y^*)\big)\,d\overline Y.
\end{array}
$$
Now we observe that
$$
|f(\overline Y)-f(\overline Y^*)|\le c|u(\overline Y)-u(Y_0)|^{\gamma}|u(\overline Y)-u(\overline Y^*)|\le ch|\overline Y-Y_0|_\sigma^{\nu\gamma},
$$
see \eqref{eq:local.holder.condition0}.
Thus
$$
|J_2|\le ch\int_{t_0}^{t_0+h}\int_{\{|\overline x-x_0|>\varrho h^{1/\sigma}\}}
\frac{|\overline Y-Y_0|_\sigma^{\nu\gamma}}{|\overline Y-Y_0|_\sigma^{N+\sigma}}\,d\overline xd\overline t\le ch^{1+\nu\gamma/\sigma}.
$$
We conclude by noting that $1+\nu\gamma/\sigma\ge\nu(1+\gamma)/\sigma$.
\qed

\begin{Lemma}\label{th:submain} Under the hypotheses of Theorem~{\rm \ref{th:main}}, $u(\cdot,t)\in
C^{\nu(1+\gamma)}(\mathbb{R}^N)$ uniformly in $t$.
\end{Lemma}

\noindent\emph{Proof. }
Once we know that  $u(x,\cdot)\in
C^{\nu(1+\gamma)/\sigma}(0,\infty)$ uniformly in $x\in\mathbb{R}^N$, we can repeat the calculations in the proof of Lemma~\ref{th:submain1}, with $\alpha$ replaced by $\nu$.
\qed

Using the worst case we can write the joint regularity in the form
\begin{equation*}
\label{eq:joint-regularity}
u\in\begin{cases}
C^{(1+\gamma)/\sigma}(Q)&\text{if } \sigma\ge1,\\[4mm]
C^{\sigma(1+\gamma)}(Q)&\text{if } \sigma\le1.
\end{cases}
\end{equation*}
with both variables playing the same role. We also have that the solution is classical since it has continuous derivatives in the sense required in the equation.

\begin{Corollary}
 Under the hypotheses of Theorem~{\rm \ref{th:main}}, the function $z:=\partial_tu=-(-\Delta)^{\sigma/2}\varphi(u)$ satisfies $z\in
C^{\nu(1+\gamma)-\sigma,(\nu(1+\gamma)-\sigma)/\sigma}_{x,t}(Q)$.
\end{Corollary}

\noindent\emph{Proof. } We point out that both sides of the equation are bounded functions and equal almost everywhere. We also know that $\partial_tu$ is H\"older continuous as a function of $t$ for a.e.~$x$, and the H\"older continuity is locally uniform. On the other hand, we easily conclude that $(-\Delta)^{\sigma/2}\varphi(u)$  is H\"older continuous as a function of $x$ for a.e.~$t$, and this happens again locally uniformly. H\"older continuity everywhere in both variables follows. \qed

Let us recall that under our assumptions $\sigma<\nu(1+\gamma)$, so that we are getting H\"older regularity for $\partial_t u$ in all cases.

\medskip

\section{Higher regularity}\label{sec.higher}
\setcounter{equation}{0}

We have already proved that solutions of~\eqref{eq:main} are  differentiable in time. However, in view of Lemma~\ref{th:submain} at this stage they are only known to be differentiable in space if $\sigma(1+\gamma)>1$, where $\gamma$ is the H\"older exponent of $\varphi'$. This assumption can we weakened.

\begin{Proposition}
\label{prop:C1} Under the assumptions of Theorem~\ref{th:main}, if $\sigma<1$ and $\gamma+\sigma>1$,  then $u\in C^{1,\alpha}(Q)$ for some $\alpha\in(0,1)$.
\end{Proposition}

\noindent\emph{Proof. } We only need to study the case $\sigma(1+\gamma)\le 1$, in which case necessarily  $\sigma<1$.

We consider the function $z=\partial_t u$, which belongs to $C^\alpha_\sigma(Q)$ for all $\alpha<\sigma$. Let $Y_0=(x_0,t_0)\in Q$ be fixed
and denote $a(Y)=\varphi'(u(Y))$, $z_0=z(Y_0)$, $a_0=a(Y_0)$. Then $z$ is a
distributional solution to the inhomogeneous fractional heat
equation
\begin{equation*}
\label{eq:inhomogeneous.fractional.he}
\partial_t z+a_0(-\Delta)^{\sigma/2}z=(-\Delta)^{\sigma/2}F_1+(-\Delta)^{\sigma/2}F_2,
\end{equation*}
where
$$
F_1=-(a-a_0)(z-z_0), \qquad F_2=-z_0a.
$$
We decompose $z$ as $z_1+z_2$,  where $z_i$ is a solution
to
$$
\partial_t z_i+a_0(-\Delta)^{\sigma/2}z_i=(-\Delta)^{\sigma/2}F_i,\qquad i=1,2.
$$
The term $z_2$ inherits the  regularity of $F_2$, that is, the regularity of $a(Y)$. As to $z_1$, we use the
fact that the function  $F_1=(a-a_0)(v-v_0)$ satisfies conditions~\eqref{eq:local.holder.condition2} and \eqref{eq:local.holder.condition}, which implies, thanks to Lemma~\ref{lem:regularity.improved2}, that $z_1$ is more smooth than $a$, hence more smooth than $z_2$. Therefore, we concentrate on the \lq bad' term, $z_2$.

The regularity of $F_2$, that is, the regularity of $\varphi'(u)$, coincides with the minimum between the regularities of $\varphi'$ and $u$. Therefore, $F_2(x,\cdot)$ belongs to $C^\gamma(\mathbb{R})$ uniformly in $x$. As for spatial regularity,  at this stage we know that $F_2(\cdot,t)$ is $C^\alpha(\mathbb{R}^N)$ uniformly in $t$ for all  $\alpha < \min\{\sigma(1+\gamma),\gamma\}$.

If $\sigma(1+\gamma)\ge\gamma$, we get $z_2(\cdot,t)\in C^{\gamma}(\mathbb{R}^N)\cap L^\infty(\mathbb{R}^N)$  uniformly in $t$, and so is $z=\partial_t u$. Using the equation we conclude that $u(\cdot,t)\in C^{\gamma+\sigma}(\mathbb{R}^N)$ uniformly in time. Since we have assumed that $\gamma+\sigma>1$, this means that $u$ is differentiable also in $x$.

If $\sigma(1+\gamma)<\gamma$,  we get $z_2(\cdot,t)\in C^{\sigma(1+\gamma)}(\mathbb{R}^N)\cap L^\infty(\mathbb{R}^N)$ uniformly in $t$. Since $w=(-\Delta)^{\sigma/2}z_2$ satisfies
$$
\partial_t w+a_0(-\Delta)^{\sigma/2}w=(-\Delta)^{\sigma/2}(-\Delta)^{\sigma/2}F_2,
$$
the regularity of $w$ is given by the regularity of $(-\Delta)^{\sigma/2}F_2$. This implies that $z_2(\cdot,t)\in C^{\alpha}(\mathbb{R}^N)$ uniformly in $t$ for all $\alpha < \min\{\sigma(2+\gamma),\gamma\}$. Repeating this argument as many times as needed, we finally obtain that $z(\cdot,t)\in C^{\alpha}(\mathbb{R}^N)$ for all $\alpha<\gamma$ uniformly in $t$. Using the equation, we obtain that $u\in C_{x,t}^{\alpha,1+\gamma}(Q)$ for all $\alpha<\sigma+\gamma$.
\qed

A similar argument allows to prove a regularity result for linear equations with variable coefficients that has an
independent interest.

\begin{Theorem}\label{th:regularity-linear} Let
$u$ be a bounded very  weak solution to $\partial_tu+(-\Delta)^{\sigma
/2}(au+b)=0$, where the coefficients satisfy $a,\,b\in
C^{1,\alpha}(Q)\cap L^\infty(Q)$, $a(x,t)\ge\delta>0$.  If $u\in
C^\alpha(Q)$ then $\partial_t u,\partial_{x_i} u\in C^\alpha (Q\cap\{t>\tau\})\cap L^\infty(Q\cap\{t>\tau\})$, $i=1,\dots,N$, for every $\tau>0$.
\end{Theorem}

\noindent\emph{Proof. }  The proof of $C^{1,\alpha}$ regularity is done by considering the linear equations satisfied by the derivatives. Boundedness for the derivatives then immediately follows, since $u\in
C^{1,\alpha}(Q\cap\{t>\tau\})\cap L^\infty(Q\cap\{t>\tau\})$.
\qed

This linear result is  used now  to obtain further regularity for the nonlinear problem, which covers in particular Theorem~\ref{th:main2}.

\begin{Theorem}
\label{cor:regularity3} If in addition to the hypotheses of Theorem~\ref{th:main}, $\varphi\in C^{k, \gamma}(\mathbb{R})$ for some $k\ge2$ and $0<\gamma<1$, then $u\in
C^{k,\alpha }(Q)$ for some $\alpha\in(0,1)$.
\end{Theorem}

\noindent\emph{Proof. }  We proceed  by induction in the derivation order. Since $\varphi\in C^{1,\theta}$ for all $\theta\in(0,1)$, Proposition~\ref{prop:C1} yields $\partial_t u,\partial_{x_i} u\in C^{\alpha}(Q)\cap L^\infty(Q)$ for some $\alpha\in(0,1)$.

Assume that the result is true for derivatives of order $j\le k-1$.
Let
$$
v_{\beta}=\partial_t^{\beta_0}\partial_{x_1}^{\beta_1}\cdots\partial_{x_N}^{\beta_N}u,\qquad \sum_{i=0}^N \beta_i=j.
$$
It is easily checked that $v_{\beta}$
satisfies an equation of the form
$$
\partial_t
v_{\beta}+(-\Delta)^{\sigma/2}(\varphi'(u)v_{\beta}+b_\beta)=0,
$$
where $b_\beta$ is a
polynomial in $v_{\beta'}$, $\beta_i'\le \beta_i$, $i=0,\dots,N$, $\sum_{i=0}^N\beta_i'\le j-1$, with coefficients involving the derivatives
$\varphi^{(l)}(u)$, $0<l\le j$. By hypothesis, $b_{\beta}\in
C^{1,\alpha}(Q)\cap L^\infty(Q)$ for some $\alpha\in(0,1)$. Since $u$ is bounded,
$a=\varphi'(u)\ge\delta>0$. Hence we may apply
Theorem~\ref{th:regularity-linear} to conclude the result.  \qed

\section{Nonlinear degenerate and singular equations}\label{sec.degenerate}
\setcounter{equation}{0}

A careful inspection of the proof of
Theorem~\ref{th:main} shows that the result has  a local
nature, and this will be exploited here to treat more general equations.

\begin{Theorem} \label{thm:local.regularity}
Let  $u$ be a bounded weak solution of
equation~\eqref{eq:main} such that $u\in C^\alpha_\sigma(\Omega)$ for some $\alpha\in(0,1)$ and some subdomain $\Omega\subset Q$. Let $\varphi\in C^{1,\gamma}(J)$, $J=(a,b)$, where $\varphi'(s)\ge c>0$ for every $s\in J$.
Under these assumptions we conclude that $\partial_t u$ and $(-\Delta)^{\sigma/2}\varphi(u)$ are H\"{o}lder continuous functions in $\mathcal{O}=\Omega\cap\{(x,t):\,u(x,t)\in J\}$. Hence $u$ is a classical solution of~\eqref{eq:main} in that set.
\end{Theorem}

\noindent\emph{Proof. }
We have to revise the proofs of all the results in Subsection \ref{sec.reg.linear} and Sections~\ref{sec.regularity.nonlinear} and \ref{sec.beyond.Holder}. For instance, in the proof of Lemma~\ref{lem:regularity.improved}, we have to replace the assumption $f\in C_\sigma^\alpha(Q)$ by $f\in C_\sigma^\alpha(\Omega)$ to conclude that $g$ belongs to the same space, and this is true since $f$ is also bounded. The same holds for Lemma~\ref{lem:regularity.improved2}.

As to Lemma~\ref{th:submain1}, we observe that the estimate of the second order differences  holds uniformly in every compact $K\subset\Omega$.
Now take a smooth cut-off function
$\phi$ with support contained in $\Omega$, with $\phi\equiv1$ in a
subset $\Omega'\subset\Omega$. Observe that the second difference of
$\psi=g\phi$ satisfies
$$
\begin{array}{l}
\psi(Y_0+Z)-2\psi(Y_0)+\psi(Y_0-Z)=\Big(g(Y_0+Z)-2g(Y_0)+g(Y_0-Z)\Big)\phi(Y_0) \\
\hspace{2cm}+\Big(\phi(Y_0+Z)-2\phi(Y_0)+\phi(Y_0-Z)\Big)g(Y_0) \\
\hspace{2cm}+\Big(\phi(Y_0)-\phi(Y_0-Z)\Big)\,\Big(g(Y_0+Z)-g(Y_0-Z)\Big)=O(h^{\alpha(1+\gamma)}),
\end{array}
$$
uniformly in $Q$, $Z=(he,0)$, $e\in\mathbb{R}^N$, $|e|=1$, $0<\alpha<\nu$. Then $(-\Delta)^{\delta/2}\psi$ is bounded for every $\delta<\alpha(1+\gamma)$, which implies that $\psi$ is $C^{\alpha(1+\gamma)}$, and thus $g\in C_x^{\alpha(1+\gamma)}(\Omega')$, for every $0<\alpha<\nu$.  The rest proceeds in the same way.
\qed

In order to apply this result we need to make sure that the solution is $C^\alpha_\sigma$ in some set $\Omega\subset Q$. This has been proved under certain conditions in~\cite{Athanasopoulos-Caffarelli}:   for some $A,\,B\in\mathbb{R}$, $A<B$, there exists a constant $C=C(A,B)>0$ such that
\begin{equation}
  \label{eq:nondeg-beta}
  \sup_{s\in[s_1,s_2]}\varphi'(s)\le C\,\frac{\varphi(s_2)-\varphi(s_1)}{s_2-s_1},\qquad A\le s_1<s_2\le
  B.
\end{equation}
Indeed, in this case every  bounded weak  solution $u$  to the equation in \eqref{eq:main} satisfying
$$
A\le\mathop{\mbox{ess inf}}_\Omega u\le \mathop{\mbox{ess
sup}}_\Omega u\le B
$$
belongs to  $C^\varepsilon(\Omega)$ for some $\varepsilon=\varepsilon(C)$, and thus $u\in C^{\alpha}_\sigma(\Omega)$, $\alpha=\nu\varepsilon$.

\medskip

\noindent\textbf{Application to the fractional porous medium
equation. } When $\varphi(u)=|u|^{m-1}u$, $m\ge1$, hypothesis
\eqref{eq:nondeg-beta} is satisfied with a constant $C=m$ which does
not depend on $A,B$. Therefore,  bounded weak solutions are
uniformly H\"older continuous in $\mathbb{R}^N\times(\tau,\infty)$,
$\tau>0$. However, the equation degenerates when $u=0$. Hence the
application of Theorem~\ref{thm:local.regularity} only yields the regularity stated there
in the set $\{u\ne0\}$.

In the fast diffusion case $m<1$ hypothesis \eqref{eq:nondeg-beta} only holds if
$A>0$ or $B<0$. Thus, $C^\alpha$ regularity is only guaranteed in the  positivity set (or
negativity set) of a solution.  Nevertheless, the application of our result
leads to the same conclusion as in the case $m>1$ in the set $\{u\ne0\}$.

On the other hand, in our paper \cite{pqrv2} we prove, for all $m>0$, that when the initial
value is nonnegative  the solution is strictly positive
everywhere for positive times. We obtain that the solution belongs
to $C^\alpha$ in this case,  and the application of the results of the present paper then imply that the solution is classical. The positivity property holds for all
$m>0$, which is in sharp contrast with the nonlinear theory with the
standard Laplacian and $m>1$, where the existence of free boundaries
is well-known~\cite{vazquez}.

\section{Theory of existence and basic properties}\label{sec.exist}
\setcounter{equation}{0}

As a complement to the previous regularity theory, we devote this
section to provide a  survey of the main facts of the existence and
uniqueness theory for the Cauchy problem for equation~\eqref{eq:main},
\begin{equation}
\label{Cauchy.problem}
\tag{CP}
\left\{
\begin{array}{ll}
\partial_tu+(-\Delta)^{\sigma/2}\varphi(u)=0  \quad&\mbox{in } Q,\\[4mm]
u(\cdot,0)=u_0 \quad&\mbox{in }\mathbb{R}^N.
\end{array}
\right.
\end{equation}
Such a theory
has been developed in great detail in the paper \cite{pqrv2} for the
case where  $\varphi$ is a power function. As in the case of the standard (local) porous medium equation,
many of the basic features of the theory can be extended to more general nonlinearities
$\varphi$, as long as they are continuous and nondecreasing, cf.~\cite{daska-kenig}. Therefore, we will outline here how such extension can be done in the fractional case $\sigma\in (0,2)$, with special attention to the points where the arguments differ.

Let us recall the concept of \emph{weak solution} to the Cauchy problem~\eqref{Cauchy.problem}: a function $u\in C([0,\infty): L^1(\mathbb{R}^N))$
such that (i) $\varphi(u) \in L^2_{\rm
loc}((0,\infty):\dot{H}^{\sigma/2}(\mathbb{R}^N))$; (ii) identity
\begin{equation*}\label{weak-nonlocal}
\displaystyle \int_0^\infty\int_{\mathbb{R}^N}u\partial_t
\zeta\,dxdt-\int_0^\infty\int_{\mathbb{R}^N}(-\Delta)^{\sigma/4}\varphi(u)(-\Delta)^{\sigma/4}\zeta\,d
xdt=0
\end{equation*}
holds for every $\zeta\in C_c^\infty(Q)$; and
(iii) $u(\cdot,0)=u_0$ almost everywhere. The (homogeneous) fractional Sobolev space $\dot{H}^{\sigma/2}(\mathbb{R}^N)$ is the space of locally integrable functions $\zeta$ such that $(-\Delta)^{\sigma/4}\zeta\in L^2(\mathbb{R}^N)$. We point that this is a convenient choice among other possible notions of weak solution, and it can be described as a weak  $L^1$ -energy solution to be specific.

\subsection{Solutions with bounded initial data} \label{sec.exist.1}

We will start by considering the theory for  initial data
$$
u_0\in L^1(\ren)\cap L^\infty(\ren).
$$
Existence and uniqueness are proved by using the definition of the  fractional Laplace operator
based in the extension technique developed by Caffarelli and Silvestre \cite{Caffarelli-Silvestre}, which is a generalization of the well-known Dirichlet to Neumann operator corresponding to $\sigma=1$.
Thus, if $g=g(x)$ is a smooth bounded
function defined in $\mathbb{R}^N$, its $\sigma$-harmonic extension
to the upper half-space, $v=\E(g)$, is  the unique smooth bounded
solution $v=v(x,y)$ to
\begin{equation}
\left\{
\begin{array}{ll}
\nabla\cdot(y^{1-\sigma}\nabla v)=0\quad &\text{in
}\mathbb{R}^{N+1}_+\equiv\{(x,y)\in\mathbb{R}^{N+1}:
x\in\mathbb{R}^N,
y>0\},\\[4mm]
v(\cdot,0)=g\quad&\text{in }\mathbb{R}^N.
\end{array}
\right. \label{sigma-extension}
\end{equation}
Then it turns out, see \cite{Caffarelli-Silvestre},
that
\begin{equation}
\label{fract-lapla}
-\dfrac{\partial v}{\partial y^\sigma}\equiv-\mu_{\sigma}\lim_{y\to0^+}y^{1-\sigma}\frac{\partial v}{\partial
y}=(-\Delta)^{\sigma/2} g(x),
\end{equation}
where $\mu_{\sigma}={2^{\sigma-1}\Gamma(\sigma/2)}/{\Gamma(1-\sigma/2)}$. In \eqref{sigma-extension} the operator  $\nabla$ acts in all
$(x,y)$ variables, while in \eqref{fract-lapla}
$(-\Delta)^{\sigma/2}$ acts only on the $x=(x_1,\cdots,x_N)$
variables.

Using this approach, problem~\eqref{Cauchy.problem} can be written in an
equivalent local form. If $u$ is a solution, then  $w=\E(\varphi(u))$
solves
\begin{equation}
\left\{
\begin{array}{ll}
\nabla\cdot(y^{1-\sigma}\nabla w)=0,\qquad &(x,y)\in\mathbb{R}^{N+1}_+,\, t>0,\\
\dfrac{\partial w}{\partial y^\sigma}-\dfrac{\partial
\beta(w)}{\partial
t}=0,\qquad&x\in\mathbb{R}^{N},\,y=0,\, t>0,\\
\beta(w)=u_0,\qquad&x\in\mathbb{R}^{N},\, y=0,\, t=0,
\end{array}
\right. \label{pp:local}
\end{equation}
where $\beta=\varphi^{-1}$.
Conversely, if we obtain a solution $w$
to~\eqref{pp:local}, then $u=\beta(w)\big|_{y=0}$ is a solution
to~\eqref{Cauchy.problem}.

We use the concept of weak solution for problem \eqref{pp:local} obtained by multiplying by a test function $\zeta$,
\begin{equation*}
\label{weak-local}
\displaystyle \int_0^\infty\int_{\mathbb{R}^{N}}\beta(w)\dfrac{\partial
\zeta}{\partial
t}\,dxdt-\mu_\sigma\int_0^\infty\int_{\mathbb{R}^{N+1}_+}y^{1-\sigma}\langle\nabla
w,\nabla \zeta\rangle\,dxdydt=0.
\end{equation*}
We then introduce the energy space
$X^\sigma(\mathbb{R}^{N+1}_+)$, the completion of
$C_c^\infty(\mathbb{R}^{N+1}_+)$ with the norm
\begin{equation*}
  \|v\|_{X^\sigma}=\left(\mu_\sigma\int_{\mathbb{R}^{N+1}_+} y^{1-\sigma}|\nabla
v|^2\,dxdy\right)^{1/2}. \label{norma2}
\end{equation*}

In order to solve the evolution problem, which is our concern, we use the Nonlinear Semigroup Generation Theorem due to Crandall-Liggett \cite{crandall-liggett}. We are thus reduced to deal with the related elliptic  problem
\begin{equation}\label{eq:local-elliptic}
\left\{
\begin{array}{ll}
\nabla\cdot(y^{1-\sigma}\nabla w)=0,\qquad &x\in\mathbb{R},\,y>0,\\
-\dfrac{\partial w}{\partial y^\sigma}+\beta(w)=g,\qquad&x\in\mathbb{R},\,y=0,
\end{array}\right.
\end{equation}
with $g\in L^1_+({\mathbb{R}})\cap L^\infty({\mathbb{R}})$. As in the case treated in \cite{pqrv2}, in order to get a solution by variational techniques, it is convenient  to replace the half space $\mathbb{R}^{N+1}_+$ by a half ball
$B_R^+=\{|x|^2+y^2<R^2,\,x\in \mathbb{R}^{N},\,y>0\}$. We impose zero Dirichlet
data on the \lq\lq new part'' of the boundary. Therefore we are led
to study the problem
\begin{equation}
\left\{
\begin{array}{ll}
\nabla\cdot(y^{1-\sigma}\nabla w)=0\quad &\mbox{in } B_R^+,\\
w=0\quad &\mbox{on }
\partial B_R^+\cap\{y>0\},\\
-\dfrac{\partial w}{\partial y^\sigma}+\beta(w)=g\quad&\mbox{on } D_R:=\{|x|<R,\, y=0\},
\end{array}
\right. \label{pp:local-elliptic-bdd}
\end{equation}
with $g\in L^\infty(D_R)$ given. Minimizing the functional
$$
J(w)=\frac{\mu_\sigma}2\int_{B_R^+}y^{1-\sigma}|\nabla
w|^2+\int_{D_R}B(w)- \int_{D_R}wg,
$$
$B'=\beta$, in the admissible set $\mathcal{A}=\{w\in H^1(B_R^+;y^{1-\sigma}): 0\le \beta(w)\le \|g\|_\infty\}$,
we obtain a unique solution $w=w_R$ to problem
\eqref{pp:local-elliptic-bdd}.
Moreover, if $g_1$ and $g_2$ are two admissible data, then the
corresponding weak solutions satisfy the $L^1$-contraction property
\begin{equation*}\label{eq:-contraction-bdd}
\int_{D_R}\left(\beta({w_1(x,0)})-\beta({w_2(x,0)})\right)_+\,dx\le
\int_{\mathbb{R}}\left(g_1(x)- g_2(x)\right)_+\,dx.
\end{equation*}

The passage to the limit $R\to\infty$ uses the
monotonicity in $R$ of the approximate solutions $w_R$. We obtain a
function $w_\infty=\lim_{R\to\infty}w_R$ which is a weak solution to
problem \eqref{eq:local-elliptic}. The above contractivity property
also holds in the limit. Moreover,
$\|\beta(w_\infty(\cdot,0))\|_{L^\infty(\mathbb{R})}\le\|g\|_{L^\infty(\mathbb{R})}$,
and  $w_\infty \ge0$, since $g\ge0$.

Now, using the Crandall-Liggett
Theorem we obtain the existence of a unique mild solution
$\overline w$ to the evolution problem \eqref{pp:local}. To prove
that $\overline w$ is moreover a weak solution to problem
\eqref{pp:local}, one needs to show that it lies in the right energy
space. This is done using the same technique as in~\cite{pqrv},
which yields the energy estimate
\begin{equation*}\label{L2grad}
\mu_\sigma\int_0^T\int_0^\infty\int_{\mathbb{R}^N} y^{1-\sigma}|\nabla \overline
w(x,y,t)|^2\,dxdydt\le \int_{\mathbb{R}^N}B(u_0(x))\,dx\quad \text{for every }T>0.
\end{equation*}
Hence the function $u=\beta({\overline
w(\cdot,0)})$ is a  weak solution to  problem~\eqref{Cauchy.problem}. In addition, $\|\beta(\overline
w(\cdot,0))\|_{L^\infty(\mathbb{R}\times(0,\infty))}\le\|u_0\|_{L^\infty(\mathbb{R})}$,
and $\overline w \ge0$.  Recalling the isometry
between  $\dot H^{\sigma/2}(\mathbb{R}^N)$ and $X^\sigma(\mathbb{R}^{N+1}_+)$, we obtain
\begin{equation*}\label{L2frac}
\int_0^T\int_{\mathbb{R}^N} |(-\Delta)^{\sigma/4}
\varphi(u)(x,t)|^2\,dxdt\le \int_{\mathbb{R}^N}B(u_0(x))\,dx\quad \text{for every }T>0.
\end{equation*}

The Semigroup Theory also guarantees that the constructed solutions satisfy the
$L^1$-contraction property $\|u(\cdot,t)-\widetilde
u(\cdot,t)\|_1\le \|u_0-\widetilde{u}_0 \|_1$.

Uniqueness follows by the standard argument due to Oleinik et al.~\cite{OKC}, using here the test function
$$
\zeta(x,t) =\left\{\begin{array}{ll}\displaystyle\int_t^T
(\varphi(u)-\varphi(\widetilde u))( x,s)\,ds,\qquad& 0\le t\le
T,\\[8pt]
0,\qquad& t\ge T,
\end{array}
\right.
$$
in the weak formulation for the difference of two solutions $u$ and $\widetilde u$.

Summarizing, we have proved  the following existence and uniqueness result.
\begin{Theorem}
Let $\varphi\in C(\mathbb{R})$ be  nondecreasing.
Given $u_0\in L^1(\ren)\cap L^\infty(\ren)$ there exists a unique bounded weak $L^1$-energy solution to problem~\eqref{Cauchy.problem}.
\end{Theorem}

\subsection{Solutions with unbounded data. Boundedness and decay}

If the (nondecreasing)
nonlinearity $\varphi$ satisfies  $\varphi'(u)\ge C |u|^{m-1}$ for
some  $m\in\mathbb{R}$  and $|u|\ge C$, then  weak solutions with
initial data in $L^1(\mathbb{R}^N)\cap L^p(\mathbb{R}^N)$, where
$p\ge1$ satisfies $p>p(m)=(1-m)N/\sigma$, become immediately bounded, hence, thanks to our results,
classical.

The idea is to take as test function  in the weak formulation $\zeta=(|u|-1)_+^{p-1}\mbox{sign}(u)$. Though  $u$
is not differentiable in time a.e.~for a general $\varphi$, this is
not needed for the proof, since a regularization procedure, using some Steklov averages, allows to
bypass this difficulty; see for example the classical
paper~\cite{Aronson-Serrin} for the case of local operators. Hence, we only have to  check that  $\zeta\in L^2_{\rm
loc}((0,\infty):\dot{H}^{\sigma/2}(\mathbb{R}^N))$ for every $p\ge2$. This will follow from the following result applied to $v=\varphi(u)$.
\begin{Proposition}
  \label{prop:powerin H-sigma}
  If $v\in\dot H^{\sigma/2}(\mathbb{R}^N)\cap L^\infty(\mathbb{R}^N)$ and $f'\in L^\infty_{\rm loc}(\mathbb{R}^N)$, then $f(v)\in\dot H^{\sigma/2}(\mathbb{R}^N)$.
\end{Proposition}

\noindent\emph{Proof. } Using the extension technique we have
$$
\begin{aligned}
  \|f(v)\|_{\dot H^{\sigma/2}(\mathbb{R}^N)}^2&=\mu_\sigma\int_{\mathbb{R}^{N+1}_+} y^{1-\sigma}|\nabla\E(f(v))|^2\,dxdy \le\mu_\sigma\int_{\mathbb{R}^{N+1}_+} y^{1-\sigma}|\nabla f(\E(v))|^2\,dxdy\\
  &\le\mu_\sigma\|f'(\E(v))\|^2_\infty\int_{\mathbb{R}^{N+1}_+} y^{1-\sigma}|\nabla \E(v)|^2\,dxdy\le c\|v\|_{\dot H^{\sigma/2}(\mathbb{R}^N)}^2.
\end{aligned}
$$
\qed

\noindent \emph{Remark. } If moreover $f$ is convex then, noting that $\|\E(v)\|_\infty\le \|v\|_\infty$, we deduce the estimate $\|f(v)\|_{\dot H^{\sigma/2}(\mathbb{R}^N)}\le\|f'(v)\|_\infty\|v\|_{\dot H^{\sigma/2}(\mathbb{R}^N)}$.

If we use the above test function and apply the generalized Stroock-Varopoulos inequality, proved in \cite[Lemma~5.2]{pqrv2}, we  obtain
$$
\int_{\ren}(|u|-1)_+^{p}(x,t_2)\,dx-\int_{\ren}(|u|-1)_+^{p}(x,\overline t)\,dx \le
-c\int_{\overline t}^{t_2}\int_{\ren} |(-\Delta)^{\sigma/4}G(u)|^2\,dxdt,
$$
where
$$
|G'(u)|^2=\varphi'(u)(|u|-1)_+^{p-2}\ge c(|u|-1)_+)^{m+p-3}.
$$
Now using the Hardy-Littlewood-Sobolev  inequality~\cite{Hardy-Littlewood},~\cite{Sobolev}  if $N>\sigma$, we get
$$
\int_{\ren}(|u|-1)_+^{p}(x,\overline t)\,dx \ge
c\int_{\overline t}^{t_2}\left(\int_{\ren} (|u|-1)_+^{\frac{(m+p-1)N}{N-\sigma}}\,dx\right)^{\frac{N-\sigma}N}dt
$$
for every $p\ge2$. 
If in addition $p>(1-m)N/\sigma$, this inequality is enough to apply a standard Moser's iteration technique  to obtain an $L^p$-$L^\infty$ smoothing effect. We can then weaken the restriction $p\ge 2$ to $p\ge 1$ using interpolation; see~\cite{pqrv2}. Take note that in the case $N=1\le\sigma<2$ we must replace the Hardy-Littlewood-Sobolev inequality by a Nash-Gagliardo-Nirenberg inequality; see \cite[Lemma~5.3]{pqrv2}. We omit further details.

Let us state precisely the smoothing result thus obtained for future reference.
\begin{Theorem}
\label{th:smoothing2}
Let $\varphi\in C^1(\mathbb{R})$ be such that
 $\varphi'(u)\ge C |u|^{m-1}$ for
some  $m\in\mathbb{R}$  and $|u|\ge C$. If $u_0\in L^1({\mathbb{R}^N})\cap
L^p({\mathbb{R}^N})$, where $p\ge1$ satisfies $p>(1-m)N/\sigma$,  then there exists a unique weak $L^1$-energy solution to the Cauchy problem~\eqref{Cauchy.problem} which is bounded in $\mathbb{R}^N\times(\tau,\infty)$ for all  $\tau>0$. This solution moreover satisfies
\begin{equation}
\sup_{x\in{\mathbb{R}^N}}|u(x,t)|\le
\max\{C, \, C_1\,t^{-\gamma_p }\|u_0\|_{p}^{\delta_p}\}
\label{eq:L-inf-p2}\end{equation}
with
$\gamma_p=N/(N(m-1)+{\sigma}p)$ and $\delta_p=\sigma p\gamma_p/N$,
the constant $C_1$ depending on $m,\,p,\,\sigma, C,$ and $N$.
\end{Theorem}

\noindent \emph{Remark. } If the function $\varphi$ satisfies the condition $\varphi'(u)\ge C |u|^{m-1}$ for every $u\in\mathbb{R}$ and some fixed $m>0$, a classical scaling argument allows to obtain a decay estimate for every $t>0$; see for instance~\cite{JLVSmoothing}. Indeed, the function $v(x,t)=\lambda^{\gamma_p}u(\lambda^{p\gamma_p/N},\lambda t)$ is a solution to the equation $\partial_tv+(-\Delta)^{\sigma/2}\widetilde\varphi(v)=0$ with $\widetilde\varphi(s)=\lambda^{m\gamma_p}\varphi(\lambda^{-\gamma_p}s)$, which satisfies the same condition on the derivative. Thus, applying~\eqref{eq:L-inf-p2} at $t=1$ and putting $\lambda=t$ we get
$$
\|u(\cdot,t)\|_\infty=t^{-\gamma_p}\|v(\cdot,1)\|_\infty\le C_1t^{-\gamma_p}\|u_0\|_p^{\delta_p}
$$
for every $t>0$.

Existence for data which are unbounded is proved by approximation; see~\cite{pqrv2} for the details in the case where the nonlinearity is a pure power. As for uniqueness, continuity in $L^1$ guarantees that two solutions with the same initial data do not differ more than $\varepsilon$ in $L^1$ norm for some small enough time. Since for positive times solutions are assumed to be bounded, we may use the $L^1$ contraction property to prove that the distance between the two solutions stays smaller than $\varepsilon$ for any later time. Since $\varepsilon$ is arbitrary, uniqueness follows.

\section{Extensions and comments}\label{sec.extension}
\setcounter{equation}{0}

\noindent\textsc{Some applications. } Equation~\eqref{eq:main} appears in the study of hydrodynamic limits of interacting particle systems with long range dynamics.  Thus, in \cite{Jara-Landim-Sethuraman}, Jara
and co-authors study the non-equilibrium functional central limit theorem for the position of a tagged
particle in a mean-zero one-dimensional zero-range process. The asymptotic behavior of the particle
is described by a stochastic differential equation governed by the solution of \eqref{eq:main}.

In several space dimensions, equations like~\eqref{eq:main} occur in boundary heat control, as already
mentioned by Athanasopoulos and Caffarelli~\cite{Athanasopoulos-Caffarelli} , where they refer to the model formulated in the book
by Duvaut and Lions~\cite{Duvaut-Lions}, and use the extension technique of Caffarelli and Silvestre.

For a more thorough discussion on applications  see~\cite{BoVa2012}.

\medskip

\noindent\textsc{Regularity for unbounded solutions. } In our proofs we are requiring the solutions to be bounded in order to make the integrals on unbounded sets convergent. However, this requirement may be not needed to this aim. It may be enough that the solutions belong to $ C([0,T]: L^1(\R^N,
\rho\,dx))$. It would be interesting to explore this possibility, since this may be helpful in the study of higher regularity.

\medskip

\noindent\textsc{Higher regularity for the fractional porous medium equation. } The main difficulty to obtain further regularity in this case is that, since the equation is not uniformly parabolic at infinity (it is not true that $0\le c\le \varphi'(u)\le C<\infty$), we do not know the derivatives to be bounded. Hence, we cannot apply Theorem~\ref{th:regularity-linear} directly. However, as mentioned in the previous paragraph, this might be circumvented by substituting the boundedness requirement by some less restrictive condition. The precise quantitative statements of the positivity property obtained in \cite{BoVa2012} might be helpful to this aim.

\medskip

\noindent\textsc{The fractional porous medium equation with sign changes. } Our results only give that the equation is satisfied in a classical sense where the solution is different from 0. It remains to determine what is the optimal regularity for changing sign solutions. A first step would be to study whether solutions are strong, i.e., whether $\partial_t u$ (and hence $(-\Delta)^{\sigma/2} u$) are functions, and not only distributions.

\medskip

\noindent\textsc{The very fast fractional porous medium equation. } The  nonlinearities
$\varphi(u)=\frac{(1+u)^{m}-1}m$, $m\ne 0$, are uniformly parabolic if we restrict ourselves to nonnegative solutions. Moreover, they fall within the hypotheses of Theorem~\ref{th:smoothing2}, if we modify the nonlinearity suitably for $u<0$, which does not matter if we only consider nonnegative solutions. Therefore, we obtain existence of $C^\infty$  solutions for all nonnegative initial data in $L^1(\mathbb{R}^N)\cap L^p(\mathbb{R}^N)$ with $p$ large enough. If
$\sigma>1-m$ and $N=1$ we can even take $p=1$.

The nonlinearity $\varphi(u)=\log(1+u)$ is also uniformly parabolic if we restrict to nonnegative solutions. In addition, after a suitable modification for $u<0$, it satisfies the hypotheses of Theorem~\ref{th:smoothing2} with $m=0$. Thus, if $N=1$ and $\sigma=1$ we are in the critical case where we need a bit more than integrability  to have existence. In~\cite{pqrv3} we proved that it is enough for $u_0$ to belong to some $L\log L$ space. The solution is then guaranteed to be $C^\infty$.

The singular nonlinearities $\varphi(u)=u^m/m$, $m<0$, and $\varphi(u)=\log u$ (with $u>0$) cannot be treated in the same way, and require new ideas.

\medskip

\noindent\textsc{The fractional Stefan problem. } For
the Stefan nonlinearity $\varphi(u)=(u-1)_+$, hypothesis
\eqref{eq:nondeg-beta} holds if $A>1$. Hence bounded weak solutions
are $C^\alpha$ in the set where $u>1$ and our main result proves that they are
$C^{1,\gamma}$, hence classical,  in that set for all $\gamma\in(0,1)$.  Let us mention that $u$ is known to be
continuous everywhere, though not $C^\alpha$. It would be interesting to determine what is the optimal
regularity for this problem.

\medskip

\noindent\textsc{Problems in fluid mechanics. } We now explore an
interesting and enlightening connection, in the case $N=\sigma=1$, between
equation~\eqref{eq:main} and Morlet's family of 1-dimensional
nonlocal transport equations with viscosity \cite{Morlet-1998},
\begin{equation}
\label{eq:v-morlet}
\partial_t v-\delta\partial_y(H(v)v)-(1-\delta)H(v)\partial_y
v=-(-\Delta)^{1/2}v,\quad 0\le\delta\le1.
\end{equation}

For a nonnegative solution $u$ to equation~\eqref{eq:main}, we consider the change of variables  $(x,t,u)\mapsto (y,\tau,v)$
given by the B\"acklund type transform
$$
y=\int_0^x(1+u(s,t))\,ds-c(t),\quad \tau=t,\qquad
v(y,\tau)=\varphi(u(x,t)),
$$
with $c'(t)=H(\varphi(u))(0,t)$. We denote $(y,\tau)=J(x,t)$. Notice
that the Jacobian of the transformation $J$ is
$\frac{\partial(y,\tau)}{\partial(x,t)}=1+u\ne0$, since $u\ge0$.
Then we may write the inverse
$$
x=\int_0^y\frac{d\theta}{1+\varphi^{-1}(v(\theta,\tau))}-\overline
c(\tau),
$$
with $\overline c\,'(\tau)=-H(\varphi(u))(0,t)/(1+u(0,t))$.

We recall that, if the operators are acting on smooth enough
functions, then the half-Laplacian $(-\Delta)^{1/2}$ can be written
in terms of the Hilbert transform
\begin{equation*} Hf(x)=\frac1\pi\mbox{
P.V.}\int_{\mathbb{R}}\frac{f(y)}{x-y}\,dy.
\end{equation*}
as
$(-\Delta)^{1/2}=H\partial_x=\partial_xH$. Therefore, we have
$$
\partial_xy=1+u,\quad \partial_ty=-H(\varphi(u))=-\widetilde H(v),
$$
where $ \widetilde H(v)=H(v\circ J)\circ J^{-1}$ is the conjugate of
the Hilbert transform $H$ by the transformation $J$. Specifically,
\begin{equation*}\label{eq:hilbert-rare}
\begin{array}{rl}
\displaystyle\widetilde
H(v(y,\tau))=H(\varphi(u(x,t)))&\displaystyle=\frac1\pi\mbox{ P.V.}
\int_{\mathbb{R}}\frac{\varphi(u(x',t))}{x-x'}\,dx'
\\ [4mm]
&\displaystyle=\frac1\pi\mbox{ P.V.}
\int_{\mathbb{R}}\frac{v(y',\tau)}{(1+\varphi^{-1}(v(y',\tau)))\int_{y'}^y\frac{d\theta}{1+\varphi^{-1}(v(\theta,\tau))}}\,dy'.
\end{array}
\end{equation*}
If $\varphi(u)=\frac{(1+u)^{m}-1}m$, $m\in[-1,0)$, then equation
\eqref{eq:main} becomes
\begin{equation*}
\partial_tv-\widetilde H(v)\partial_yv=-(1+mv)\partial_y\widetilde H(v).
\label{eq:morlet}
\end{equation*}
If instead of $\widetilde H$ we had the standard Hilbert transform
$H$, and we take $m=-\delta$, we would have an equation in Morlet's
family \eqref{eq:v-morlet}. The connection also works for the  case
$m=0$, if we take $\varphi(u)=\log(1+u)$; see~\cite{pqrv3}.

 \vskip 1cm


\noindent {\large\bf Acknowledgments}

\noindent FQ, AR, and JLV  partially supported by the Spanish project MTM2011-24696. AdP partially supported by
the Spanish project MTM2011-25287.

\vskip 1cm




\

\noindent{\bf Addresses:}

\noindent{\sc J.~L. V\'{a}zquez: }
Departamento de Matem\'{a}ticas, Universidad Aut\'{o}noma de Madrid, 28049
Madrid, Spain. (e-mail: juanluis.vazquez@uam.es).

\noindent{\sc A. de Pablo: } Departamento de Matem\'{a}ticas, Universidad Carlos III de Madrid, 28911 Legan\'{e}s,
Spain. (e-mail: arturo.depablo@uc3m.es).

\noindent{\sc F. Quir\'{o}s: }
Departamento de Matem\'{a}ticas, Universidad Aut\'{o}noma de Madrid, 28049 Madrid, Spain.
(e-mail: fernando.quiros@uam.es).

\noindent{\sc A. Rodr\'{\i}guez: }
Departamento de Matem\'{a}tica, ETS Arquitectura, Universidad Polit\'{e}cnica de Madrid, 28040 Madrid, Spain. (e-mail: ana.rodriguez@upm.es).

\end{document}